\documentclass[12 pt]{amsart}
\usepackage{amssymb, amsmath, amsfonts, amsthm, graphics,mathrsfs, mathtools}
\usepackage[usenames, dvipsnames]{xcolor}
\usepackage[hmargin=1 in, vmargin = 1 in]{geometry}
\usepackage{tikz}
\usetikzlibrary{matrix, calc, arrows, positioning, backgrounds} 
\usepackage{hyperref}
\usepackage[capitalise,noabbrev]{cleveref}
\usepackage{ytableau}

\definecolor{darkblue}{rgb}{0.0,0,0.7}

\newcommand{\SL}{\mathrm{SL}}

\newcommand{\SYT}{\mathrm{SYT}}
\newcommand{\latset}{\mathcal{P}}

\newcommand{\N}{\mathsf{N}}
\renewcommand{\S}{\mathsf{S}}
\newcommand{\E}{\mathsf{E}}
\newcommand{\W}{\mathsf{W}}

\newtheorem{Theorem}{Theorem}[section]

\newtheorem{definition}[Theorem]{Definition}
\newtheorem{proposition}[Theorem]{Proposition}

\newtheorem{theorem}[Theorem]{Theorem}

\numberwithin{equation}{section}

\theoremstyle{remark}

\title{An exceptional equinumerosity of lattice paths and Young tableaux}
\author[Ayres]{Liam Ayres}
\address[LA]{Department of Pure Mathematics, University of Waterloo, Waterloo ON N2L 3G1}
\email{liam.ayres@uwaterloo.ca}

\author[Bialo]{Evan Bialo}
\address[EB]{Department of Statistics \& Actuarial Science, University of Waterloo, Waterloo ON N2L 3G1}
\email{ebialo@uwaterloo.ca}

\author[Cook]{Aidan Cook}
\address[ACook]{Department of Pure Mathematics, University of Waterloo, Waterloo ON N2L 3G1}
\email{aidan.cook@uwaterloo.ca}

\author[Chen]{Alwin Chen}
\address[AChen]{Department of Pure Mathematics, University of Waterloo, Waterloo ON N2L 3G1}
\email{z786chen@uwaterloo.ca}

\author[Froese]{Matteus Froese}
\address[MF]{Department of Combinatorics \& Optimization, University of Waterloo, Waterloo ON N2L 3G1}
\email{mjfroese@uwaterloo.ca}

\author[Liu]{Erica Liu}
\address[EL]{Department of Combinatorics \& Optimization, University of Waterloo, Waterloo ON N2L 3G1}
\email{e237liu@uwaterloo.ca}

\author[Mohammadi Yekta]{Maryam Mohammadi Yekta}
\address[MMY]{Department of Combinatorics \& Optimization, University of Waterloo, Waterloo ON N2L 3G1}
\email{m5mohammadiyekta@uwaterloo.ca}

\author[Pechenik]{Oliver Pechenik}
\address[OP]{Department of Combinatorics \& Optimization, University of Waterloo, Waterloo ON N2L 3G1}
\email{oliver.pechenik@uwaterloo.ca}

\author[Wong]{Benjamin Wong}
\address[BW]{Department of Combinatorics \& Optimization, University of Waterloo, Waterloo ON N2L 3G1}
\email{b62wong@uwaterloo.ca}

\date{\today}
\keywords{}

\makeatletter
\@namedef{subjclassname@2020}{%
  \textup{2020} Mathematics Subject Classification}
\makeatother

\subjclass[2020]{05A19; 05A15}

\AtBeginDocument{%
   \def\MR#1{}
}

\begin{document}

\begin{abstract}
We consider families $\latset_n$ of plane lattice paths enumerated by Guy, Krattenthaler, and Sagan (1992). We show by explicit bijection that these families are equinumerous with the set $\SYT(n+2,2,1^n)$ of standard Young tableaux.
\end{abstract}

\maketitle

\section{Introduction}\label{sec:intro}
\ytableausetup{boxsize=1em}

We consider the set $\latset_n$ of lattice paths from $(0,0)$ to $(n,n)$ that
\begin{itemize}
	\item only use the steps $\N = (0,1), \S = (0,-1), \E = (1,0),$ and $\W = (-1,0)$;
	\item stay weakly inside the first quadrant of the plane;
	\item have length $2n+2$.
\end{itemize} 
It was shown by Richard K.\ Guy, Christian Krattenthaler, and Bruce E. Sagan
\cite{Guy.Krattenthaler.Sagan} that 
\begin{equation}\label{eq:enumer}
|\latset_n| = \binom{2n}{n}\frac{(4n+4)(2n+1)}{n+2}.
\end{equation}
This integer sequence appears as \cite{oeis_seq}. The same enumerations also appear in \cite{Connor.Fewster} in the context of integrals of incomplete beta functions.

What appears not to have been previously noticed is that \eqref{eq:enumer} is also the enumeration of standard Young tableaux of shape $\theta^{(n)} = (n+2,2,1^{n})$, an easy consequence of the hook-length formula \cite{Frame.Robinson.Thrall} for counting standard Young tableaux of any fixed shape. (The shape $\theta^{(n)}$ is an instance of the \emph{near-hook shapes} studied in \cite{Langley.Remmel}.) 
That is, we have 
\begin{equation}\label{eq:enumer_tabs}
\left| \SYT \left(\theta^{(n)} \right) \right| = \binom{2n}{n}\frac{(4n+4)(2n+1)}{n+2}.
\end{equation}
Here we give a direct proof of this equinumerosity by exhibiting an explicit bijection between the sets $\latset_n$ of lattice paths and $\SYT(\theta^{(n)})$ of standard tableaux. 

Special enumerations of families of standard Young tableaux are often a sign of deeper algebraic structure, which can often be revealed by bijecting tableaux of these families with other combinatorial objects. A famous example is the set $\SYT(k,k)$ of $2$-row rectangular tableaux, which is enumerated by the Catalan numbers; an equivariant bijection to noncrossing matchings yields a \emph{cyclic sieving} formula for the orbit structure of Sch\"utzenberger's \cite{Schutzenberger} \emph{promotion} operator on $\SYT(k,k)$ (see \cite{Reiner.Stanton.White,Petersen.Pylyavskyy.Rhoades,Rhoades} for discussion). Similarly, the orbit structure of promotion on $\SYT(k,k,k)$ and $\SYT(k,k,k,k)$ may be understood via exceptional bijections to Kuperberg's \emph{$\SL_3$-webs} \cite{Kuperberg, Khovanov.Kuperberg,Petersen.Pylyavskyy.Rhoades} and \emph{$4$-hourglass plabic graphs} \cite{Gaetz.Pechenik.Pfannerer.Striker.Swanson}. 

Another class of standard tableaux with a special enumeration is $\SYT(n-k,n-k,1^k)$, which was shown by Stanley \cite{Stanley:bij} to biject with polygon dissections of an $(n+2)$-gon by $n-k-1$ diagonals. Here, the bijection does not explain the promotion action on the set $\SYT(n-k,n-k,1^k)$. Instead, there are further bijections to rectangular \emph{increasing} tableaux and to noncrossing matchings without singleton blocks and it is the $K$-theoretic promotion on these increasing tableaux that is explicated through these bijections \cite{Pechenik:CSP}. These perspectives lead to a surprising diagrammatic basis for the Specht module for the partition shape $(n-k,n-k,1^k)$ \cite{Rhoades:skein,Patrias.Pechenik.Striker,Kim.Rhoades,Fraser.Patrias.Pechenik.Striker}.

We do not know such an application to representation theory of the bijection for $\SYT(\theta^{(n)})$ given here, but in light of this bijection and the above results, we would suggest that the Specht module for shape $\theta^{(n)}$ receive further combinatorial study. An interesting feature of the correspondence between $\SYT(\theta^{(n)})$ and $\latset_n$ is that the number of entries in each tableau does not match the number of steps in the corresponding lattice path; the lattice paths each have length $2n+2$, while the tableaux each have $2n+4$ boxes. This mismatch makes the correspondance rather subtle.

\section{Notations}\label{sec:notations}
To fix notation and conventions, we recall some standard notions in tableau combinatorics. For further background, see, for example, the textbooks \cite{Fulton:YT,Stanley:EC2}. The only non-standard content in this section is \cref{def:body_parts}.

We write $\theta^{(n)}$ for the integer partition $(n+2,2,1^n)$ of the number $2n+4$ into $n+2$ parts. We conflate $\theta^{(n)}$ with its \emph{Young diagram} in English orientation, so that the row of length $n+2$ is at the top. For example, we draw $\theta^{(3)}$ as $\ydiagram{5,2,1,1,1}$. 

\begin{definition}\label{def:body_parts}
	We refer to the first row of $\theta^{(n)}$ as its \emph{arm}, the first column as its \emph{leg}, and the unique box that is in neither its first row nor its first column as its \emph{heart}. 
\end{definition}

A \emph{standard Young tableau} $T \in \SYT(\theta^{(n)})$ is a bijective filling of the boxes of the Young diagram with the integers $1, \dots, 2n+4$ such that the entries increase along rows left to right and increase down columns top to bottom. For example, $\ytableaushort{1245,37,6,8} \in \SYT(\theta^{(2)})$. Here, the entries of the arm are $1,2,4,5$, the entries of the leg are $1,3,6,8$, and the heart entry is $7$.

For conciseness, we write our lattice paths $P \in \latset_n$ as words in the alphabet $\{\N, \S, \E, \W \}$ instead of drawing the paths in the plane. For example, $\N \N \E \S \in \latset_1$ denotes the lattice path that looks like
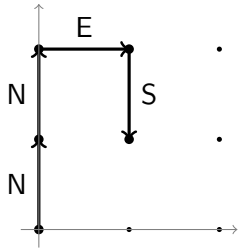
\[
\begin{tikzpicture}[scale=1.2, thick]
  \coordinate (A) at (0,0);
  \coordinate (B) at (0,1);
  \coordinate (C) at (0,2);
  \coordinate (D) at (1,2);
  \coordinate (E) at (1,1);
  \foreach \pt in {A,B,C,D,E}
    \fill (\pt) circle (1.5pt);
  \fill (1,0) circle (0.8pt);
  \fill (2,0) circle (0.8pt);
  \fill (2,1) circle (0.8pt);
  \fill (2,2) circle (0.8pt);
  \draw[very thick,->] (A) -- (B) node[midway,left] {\small $\mathsf N$};
  \draw[very thick,->] (B) -- (C) node[midway,left] {\small $\mathsf N$};
  \draw[very thick,->] (C) -- (D) node[midway,above] {\small $\mathsf E$};
  \draw[very thick,->] (D) -- (E) node[midway,right] {\small $\mathsf S$};
  \draw[->, thin, gray] (-0.2,0) -- (2.2,0);
  \draw[->, thin, gray] (0,-0.2) -- (0,2.5);
\end{tikzpicture}
\]
We refer to the steps $\N$ and $\E$ as \emph{forward} steps and refer to the steps $\S$ and $\W$ as \emph{backwards} steps. Note that the length and endpoint conditions for $\latset_n$ force that each $P \in \latset_n$ must have exactly one backwards step.

\section{First bijection}\label{sec:first_bij}

In this section, we describe and prove the correctness of an explicit bijection \[\psi : \SYT(\theta^{(n)}) \to \latset_n.\]

Let $T \in \SYT(\theta^{(n)})$. We construct the corresponding path $\psi(T) = P = p_1 p_2 \ldots p_{2n+2} \in \latset_n$ by noting the positions of the entries of $T$ in increasing order.
	The entry $1$ appears in the same position in all tableaux; hence, it carries no information and we ignore it entirely.
	
As mentioned at the end of \Cref{sec:notations}, a path in $\latset_n$ must have exactly one backwards step, either $\W$ or $\S$. We construct the path $P$ to explicitly have exactly one backwards step. If the value $2$ appears in the arm of $T$, the backwards step in $P$ will be $\S$. If instead $2$ appears in the leg of $T$, the backwards step will be $\W$.
	
	Now, for each $3 \leq i \leq 2n+4$, we construct a step of the lattice path $P$. If $i$ lies in the arm, then set $p_{i-2} = \E$. If $i$ lies in the leg, then set $p_{i-2} = \N$. If $i$ is in the heart, then $p_{i-2} \in \{\S, \W \}$ is the backwards step, with the type of this backwards step determined earlier by observing the position of the entry $2$.
	
\begin{figure}[htbp]
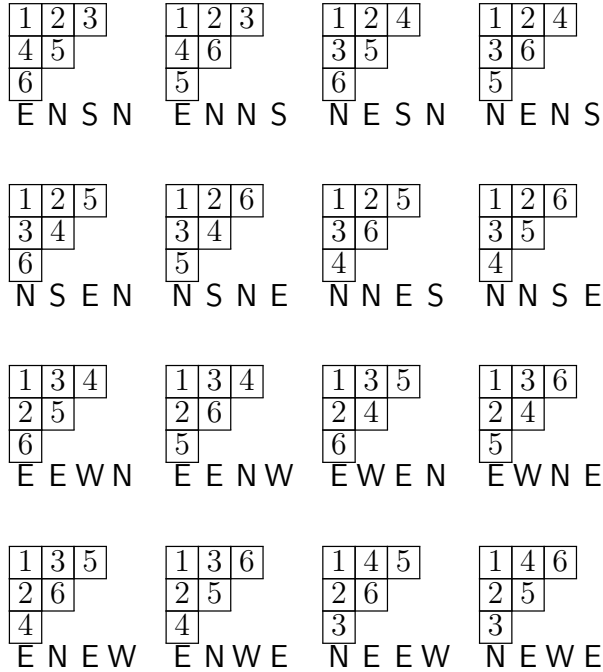

	\[\def\arraystretch{3.2}
	\begin{array}{cccc}
		\ytableaushort{123,45,6,{\none[\E]}{\none[\N]}{\none[\S]}{\none[\N]}} & \ytableaushort{123,46,5,{\none[\E]}{\none[\N]}{\none[\N]}{\none[\S]}} & \ytableaushort{124,35,6,{\none[\N]}{\none[\E]}{\none[\S]}{\none[\N]}} & \ytableaushort{124,36,5,{\none[\N]}{\none[\E]}{\none[\N]}{\none[\S]}} \\[5pt]
		\ytableaushort{125,34,6,{\none[\N]}{\none[\S]}{\none[\E]}{\none[\N]}} & \ytableaushort{126,34,5,{\none[\N]}{\none[\S]}{\none[\N]}{\none[\E]}} & \ytableaushort{125,36,4,{\none[\N]}{\none[\N]}{\none[\E]}{\none[\S]}} & \ytableaushort{126,35,4,{\none[\N]}{\none[\N]}{\none[\S]}{\none[\E]}} \\[5pt]
		\ytableaushort{134,25,6,{\none[\E]}{\none[\E]}{\none[\W]}{\none[\N]}} & \ytableaushort{134,26,5,{\none[\E]}{\none[\E]}{\none[\N]}{\none[\W]}} & \ytableaushort{135,24,6,{\none[\E]}{\none[\W]}{\none[\E]}{\none[\N]}} & \ytableaushort{136,24,5,{\none[\E]}{\none[\W]}{\none[\N]}{\none[\E]}} \\[5pt]
		\ytableaushort{135,26,4,{\none[\E]}{\none[\N]}{\none[\E]}{\none[\W]}} & \ytableaushort{136,25,4,{\none[\E]}{\none[\N]}{\none[\W]}{\none[\E]}} & \ytableaushort{145,26,3,{\none[\N]}{\none[\E]}{\none[\E]}{\none[\W]}} & \ytableaushort{146,25,3,{\none[\N]}{\none[\E]}{\none[\W]}{\none[\E]}}
	\end{array}
	\]
		\caption{The $16$ standard Young tableaux in $\SYT(\theta^{(1)})$. Below each tableau $T$ appears the corresponding lattice path $\psi(T) \in \latset_1$.}\label{fig:psi_ex}
	\end{figure}

	In summary, $\psi(T) = p_1 p_2 \ldots p_{2n+2}$ is given by
	\begin{equation}\label{eq:psi}
	p_i = \begin{cases}
		\E, & \text{if $i+2 \in \mathrm{arm}(T)$;}\\
		\N, & \text{if $i+2 \in \mathrm{leg}(T)$;}\\
		\S, & \text{if $i+2 \in \mathrm{heart}(T)$ and $2\in \mathrm{arm}(T)$;}\\
		\W, & \text{if $i+2 \in \mathrm{heart}(T)$ and $2\in \mathrm{leg}(T)$.}\\
	\end{cases}
	\end{equation}
	For an example of the map $\psi$ applied to all $16$ tableaux in $\SYT(\theta^{1})$, see \cref{fig:psi_ex}.
	
	\begin{theorem}\label{thm:first_bij}
		The map $\psi : \SYT(\theta^{(n)}) \to \latset_n$ given in \eqref{eq:psi} is a bijection.
	\end{theorem}
	
	Before proceeding with the proof of \cref{thm:first_bij}, we note that the definition of $\psi$ is rather delicate. For example, the reader may enjoy verifying that swapping the conditions for the steps $\S$ and $\W$ would generate lattice paths that leave the first quadrant of the plane and hence are not in $\latset_n$.
	
\begin{proof}[Proof of \cref{thm:first_bij}]
	The primary aspect that needs proof is the well-definedness of the map $\psi$. Let $T \in \SYT(\theta^{(n)})$. Note that a word in $\{\N, \S, \E, \W \}$ of length $2n+2$ giving a lattice path from $(0,0)$ to $(n,n)$ must contain exactly one instance of either $\S$ or $\W$; moreover, if it contains $\S$, then the other letters must be $n+1$ copies of $\N$ and $n$ copies of $\E$, while if it contains $\W$, then the other letters must be $n$ copies of $\N$ and $n+1$ copies of $\E$.
	
	By construction, $\psi(T)$ has exactly one instance of $\S$ if $2$ is in the arm and exactly one instance of $\W$ if $2$ is in the leg. (Note that by the increasingness conditions on tableaux that $2$ cannot appear in the heart nor in the intersection of the leg and arm.) In the case where $2$ is in the arm of $T$ so that $\S$ is in $\psi(T)$, there are $n$ numbers from $\{3,  \dots, 2n + 4\}$ in the arm and $n+1$ numbers from $\{3,  \dots, 2n + 4\}$ in the leg, so that $\psi(T)$ has $n$ copies of $\E$ and $n+1$ copies of $\N$, as desired. In the case where $2$ is in the leg of $T$ so that $\W$ is in $\psi(T)$, there are $n$ numbers from $\{3,  \dots, 2n + 4\}$ in the leg and $n+1$ numbers from $\{3,  \dots, 2n + 4\}$ in the arm, so that $\psi(T)$ has $n+1$ copies of $\E$ and $n$ copies of $\N$, again as desired. Thus, $\psi(T)$ gives a lattice path of the correct length, ending at the correct position.
	
	We now verify that these lattice paths remain inside the first quadrant. We must show that, if $\S$ appears, it appears after an instance of $\N$, while, if $\W$ appears, it appears after an instance of $\E$. Suppose $p_i = \S$. Then $i + 2 \in \mathrm{heart}(T)$ and $2 \in \mathrm{arm}(T)$. By increasingness of the second row of $T$, the entry $j$ appearing directly left of the $i+2$ in the heart must satisfy $j < i+2$. Moreover,  we know the locations of $1$ and $2$ in $T$, so $j \geq 3$. It follows that $p_{j-2} = \N$ is a letter of $\psi(T)$ appearing before $p_i$, as required. Similarly if $p_i = \W$, then $i + 2 \in \mathrm{heart}(T)$ and $2 \in \mathrm{leg}(T)$, so that the entry $k$ directly above the $i+2$ in the heart satisfies $3 \leq k < i + 2$. We conclude that $p_{k-2} = \E$ is a letter of $\psi(T)$ appearing before $p_i$. This completes the proof that $\psi$ is well-defined.
	
	The injectivity of $\psi$ is clear. Surjectivity follows from the existence of an inverse map $\phi$. Given $P = p_1 \ldots p_{2n+2} \in \latset_n$, we produce $T = \phi(P)$ as follows. Place value $1$ in the intersection of the arm and leg of $\theta^{(n)}$. The elements of $\{ i+ 2 : p_i = \E \}$ go in the arm of $T$, while the elements of $\{ i+ 2 : p_i = \N \}$ go in the leg of $T$. For $p_j$ the unique backwards step of $P$, place $j+2$ in the heart of $T$. If $p_j = \S$, include $2$ in the arm of $T$, while if $p_j = \W$, include $2$ in the leg of $T$. It is straightforward to check that this produces a valid standard Young tableau of the desired shape and that $\phi, \psi$ are mutually inverse.
 \end{proof}

\section{Additional bijections}

Given two sets of cardinality $k$, there are $k!$ distinct bijections between them. However, for sets of combinatorial interest, there is often one (or perhaps a small number) of these bijections that are understood to be the ``best'' or ``most natural'' ones. For example, one might desire the bijection to preserve important weight functions on the sets or to preserve the action of some group. While we find the bijection $\psi$ of \cref{sec:first_bij} attractive, we are not yet certain what properties one most wants a bijection $\SYT(\theta^{(n)}) \xrightarrow{\sim} \latset_n$ to preserve. In this section, we describe without proof several additional such bijections for possible future application.

Given a bijection $\kappa : \SYT(\theta^{(n)}) \xrightarrow{\sim} \latset_n$, one may obtain another $\kappa^\top$ by first transposing each tableau $T$ by reflecting across the main diagonal and then applying $\kappa$. These bijections $\kappa, \kappa^\top$ are essentially the same up to convention choices. We now proceed to give another bijection that, while similar in flavour to the bijection $\psi$, is fundamentally distinct.

\subsection{Second bijection}
Define a map $\xi : \SYT(\theta^{(n)}) \to \latset_n$ as follows. Let $T \in \SYT(\theta^{(n)})$ and let $H$ be the entry in the heart of $T$. We define $\xi(T) = p_1 p_2 \ldots p_{2n+2}$ by
\begin{equation}\label{eq:xi}
	p_i = \begin{cases}
		\E, & \text{if $i+1 \in \mathrm{arm}(T)$ and $i+2 < H$}; \\ 
		\N, & \text{if $i+1 \in \mathrm{leg}(T)$ and $i+2 < H$}; \\ 
		\S, & \text{if $i+1 \in \mathrm{arm}(T)$ and $i+2 = H$}; \\ 
		\W, & \text{if $i+1 \in \mathrm{leg}(T)$ and $i+2 = H$}; \\ 
		\E, & \text{if $i+2 \in \mathrm{arm}(T)$ and $H  < i + 2$}; \\ 
		\N, & \text{if $i+2 \in \mathrm{leg}(T)$ and $H < i+2$}. \\ 
	\end{cases} 
\end{equation}

Then we have the following.

\begin{proposition}
	The map $\xi : \SYT(\theta^{(n)}) \to \latset_n$ given in \eqref{eq:xi} is a bijection. \qed 
\end{proposition}

Note that the backwards step of $\psi(T)$ is determined by the entry of the heart of $T$ and the position of $2$, whereas the backwards step of $\xi(T)$ is determined by the entry of the heart of $T$ and the numerically previous entry. Of course, $\xi^\top$ yields yet another bijection.

\section*{Acknowledgements}
The paper originates in a course CO 431/631 ``Symmetric Functions'' taught by Pechenik in Fall 2024 at the University of Waterloo. The other authors were all students in this course.

Liu was partially supported by Discovery Grant RGPIN-2021-02382 from the Natural Sciences and Engineering Research Council (NSERC) of Canada.
Mohammadi Yekta was partially supported by Jonathan Leake's NSERC Discovery Grant.
Pechenik was partially supported by an NSERC Discovery Grant (RGPIN-2021-02391) and Launch Supplement (DGECR-2021-00010). Wong was partially supported by an NSERC Alexander Graham Bell CGSM as well as David Gosset’s NSERC Discovery Grant.

\bibliographystyle{amsalpha} 
\bibliography{431bij}
\end{document}